\documentclass[a4paper,abstracton]{scrartcl}
\usepackage[latin1]{inputenc}
\usepackage{amsfonts,amsmath,amsthm,amssymb,eucal,url}
\usepackage[dvips,ps,matrix,arrow,curve]{xy}
\usepackage[english]{babel}
\usepackage[T1]{fontenc}

\CompileMatrices

\newtheorem{theor}{Theorem}

\newtheorem{prop}[theor]{Proposition} 
\newtheorem{lemma}[theor]{Lemma}

\theoremstyle{definition}
\newtheorem*{df}{Definition}

\theoremstyle{remark}
\newtheorem{rem}[theor]{Remark}

\newcommand{\abs}[1]{\left\lvert#1\right\rvert} % |x|

\newcommand{\field}[1]{\mathbf{#1}}
\newcommand{\Q}{\field{Q}}
\newcommand{\Z}{\field{Z}}
\newcommand{\R}{\field{R}}
\newcommand{\C}{\field{C}}

\newcommand{\N}{\field{N}}

\newcommand{\bigoh}[1]{O\left(#1\right)} % O() notation
\newcommand{\smalloh}[1]{o\left(#1\right)} % o() notation
\renewcommand{\Re}{\operatorname{Re}}
\renewcommand{\phi}{\varphi}
\renewcommand{\epsilon}{\varepsilon}
\renewcommand{\leq}{\leqslant}
\renewcommand{\geq}{\geqslant}

\titlehead{Dipartimento di Matematica,
  Università di Milano\\
  Via Saldini 50,
  20130 Milano, Italia\\
  e-mail: \url{Ottavio.Rizzo@mat.unimi.it}}

\title{The height measure of $p$-adic balls}
\author{Ottavio G. Rizzo}
\begin{document}

\maketitle

\begin{abstract}
  \noindent In this paper we give the height measure of $p$-adic balls. In
  other words, given any $x\in\Q_p$, we give the chance that a random
  rational number $r$ satisfies $\abs{r-x}_p \leq \epsilon$. \smallskip

  \noindent{\sectfont MSC 2000:}\enspace 11B05, 11G50, 11S80, 28C10
\end{abstract}

Let $H(m/n) = \max\{\abs{m},n\}$ be the \emph{height} of the rational
number $m/n$, where $m\in\Z$, $n\in\N$ and $\gcd(m,n)=1$. Given
$U\subseteq\Q$, consider the limit
\begin{equation*}
  \lim_{t\to\infty} 
  \frac{\#\{ r\in U : H(r)\leq t\}}{\#\{r\in\Q : H(r)\leq t}.
\end{equation*}
If it exists, we denote its value $\mu(U)$: the \emph{height density} of
$U$.

More in general, if $U$ is a subset of a completion $K$ of $\Q$ and the limit
\begin{equation}
  \label{eq:def}
  \mu(U) = \lim_{t\to\infty} \frac{%
    \# \{r\in \bar U \cap \Q : H(r) \le t \}}{%
    \# \{r\in \Q : H(r) \le t \}}
\end{equation}
exists, we say that $U$ is $\mu$-measurable and call $\mu(U)$ the
\emph{height measure} of $U$.

In \cite{rizzo99:_averag} we proved that any interval $(a,b)\subset\R$ is
$\mu$-measurable and gave a simple formula for its value. If $p$ is a
finite place of $\Q$, denote as usual $v_p$ the associated valuation,
$\abs{\cdot}_p$ the norm and $\Q_p$ the completion. Let $B(x,p^{-e}) =
\{r\in\Q_p : \abs{x-r}_p \leq p^{-e}\}$ be the (closed) $p$-adic
\emph{ball} of centre $x$ and radius $p^{-e}$. In this paper we prove that:
\begin{theor}
  \label{theor:1}
  For any $x\in\Q_p$, the $p$-adic ball $B(x,p^{-e})$ is
  $\mu$-measurable. Moreover: 
  \begin{itemize}
  \item if $e\leq v(x)$:
    \begin{equation*}
      \mu\bigl( B(x,p^{-e}) \bigr) =
      \begin{cases}
        \displaystyle \frac{p^{1-e}}{p+1_{\phantom {\bigl |}}}
        & \text{if $e\geq0$,}\\
        \displaystyle 1 - \frac{p^{e}}{p+1} & \text{if $e\leq0$;}
      \end{cases}
    \end{equation*}
  \item while, if $e>v(x)$:
    \begin{equation*}
      \mu\bigl( B(x,p^{-e}) \bigr) =
      \begin{cases}
        \displaystyle \frac{p^{1-e}}{p+1_{\phantom {\bigl |}}} 
        & \text{if $v(x)\geq0$} \\
        \displaystyle
        \frac{p^{1-e}}{\abs{x}_p^2(p+1)} 
        & \text{if $v(x)<0$} \\
      \end{cases}
    \end{equation*}
  \end{itemize}
  In particular,
  \begin{equation*}
    \mu(\{r\in\Q : v(r-x) = e \}) =
    \begin{cases}
      \displaystyle
      p^{-\abs{e}} \frac{p-1}{p+1} & \text{if $e\leq v(x)$ or $e\geq0$;} \\
      \displaystyle
      \frac{p}{p+1}\left(
        \frac{p}{\abs{x}_p^2} - 1
      \right) & \text{if $v(x)< e = -1$;} \\
      \displaystyle
      \frac{p^{-e}}{\abs{x}_p^2}\, \frac{p-1}{p+1} & \text{if $v(x)< e < -1$.}
    \end{cases}
  \end{equation*}
\end{theor}

\section{A dutiful note on measure theory}
\label{sec:some-measure-theory}

Unfortunately, Eq.~\eqref{eq:def} does not define what is usually called a
measure: indeed, $\Q$ is not $\sigma$-additive, so no function on it may be
$\sigma$-additive. For example,
\begin{math}
  \Q = \bigcup_{r\in \Q} \{r\},
\end{math}
but 
\begin{equation*}
  1 = \mu(\Q) \neq \sum_{r\in\Q} \mu(\{r\}) = 0.
\end{equation*}

This is not a serious problem, since we can easily define a real (pun
intended!) measure on $\Q_p$ which agrees, on $p$-adic balls, with our
definition: 

\begin{df}
  Let $p$ be a finite or infinite place of $\Q$.
  For any $E \subset \Q_p$ and $\delta>0$, let
  \begin{equation*}
    \mu_\delta(E) = \inf_{\substack{
        \abs{B_i} \leq \delta \\
        \bigcup B_i \supset E
        }} \sum \mu(B_i),
  \end{equation*}
  where the $B_i$ are $p$-balls and the unions are
  countable. Furthermore, let
  \begin{equation*}
    \mu^*(E) = \sup_{\delta>0} \mu_\delta(E).
  \end{equation*}
\end{df}

\begin{theor}
  \label{thm:3}
  For any place $p$ of $\Q$,
  the set function $\mu^*$ is a $\sigma$-additive measure on $\Q_p$,
  the Borel sets are measurable and 
  \begin{math}
    \mu^*(B) = \mu(B)
  \end{math}
  for any $p$-ball.
\end{theor}
\begin{proof}
  Recall that $\mu$ is an additive set function by theorem 4 of
  \cite{rizzo99:_averag}, hence by general measure theory on metric
  spaces (see for example theorem 23 of \cite{rogers70:_hausd_measur})
  $\mu^*$ is a $\sigma$-additive measure on $\Q_p$ and the Borel sets
  are measurable.
  
  We are left to prove that $\mu^*$ coincides with $\mu$ on $p$-balls.
  If $p=\infty$, since by theorem 4 of \cite{rizzo99:_averag} $\mu(B)$
  is essentially the length of the interval $B$, this is a classical
  result: see, for example, §5 of \cite{rogers70:_hausd_measur}.
  Suppose now that $p<\infty$: we claim that
  \begin{equation}
    \label{eq:1}
    \mu_\delta(B) = \mu(B), \qquad \text{for every $\delta>0$.}
  \end{equation}
  Fix such a $\delta$, and let $\{B_i\}_{i\in I}$ be $p$-balls with
  $\abs{B_i} \leq \delta$ and $\bigcup B_i \supset B$. Since $B$ is
  compact and each $B_i$ is open, we may suppose that $I$ is finite.
  For every $i\in I$, we may clearly suppose that $B_i \cap B \neq
  \emptyset$; if $x_i \in B_i \cap B$, then we can take $x_i$ to be
  the centre of both; so either $B_i \subset B$ or $B_i \supset B$,
  the latter being an uninteresting case. Similarly, we may assume
  that the $B_i$ are all pairwise disjoint. Therefore, we have
  \begin{math}
    B = \bigcup B_i
  \end{math}
  where the union is finite and disjoint; equation~\eqref{eq:1} now
  follows from the additivity of $\mu$.
\end{proof}

\section{Preliminary results}
\label{sec:preliminary-results}

From now on we fix a finite prime $p$.

\begin{df}
  In analogy to Euler's $\phi$ function, we define for any positive integer
  $t$ and any positive number $x$, a function
  \begin{equation*}
    \phi(t,x) = \# \{\text{positive integers $\leq x$ which are relatively
    prime to $t$}\}
  \end{equation*}
\end{df}
We proved in \cite{rizzo99:_averag} that:
\begin{prop}
  \label{p:3}
  Denote $d(n)$ the number of divisors of $n$.
  Then, for any $x, t>0$ we have that
  \begin{equation*}
    \phi(t,x) = \frac{x}{t}\phi(t) \pm d(t),
  \end{equation*}
  where $a=b\pm\delta$ means that $\abs{a-b}\leq\delta$.
\end{prop}

\begin{lemma}
  \label{l:1}
  Suppose $t$, $a$, $e$ are integer numbers with $t > 1$
  and fix $T>0$. Then
  \begin{multline*}
    \# \{ n\in\Z : 0 \leq n \leq T, \gcd(n,t)=1, \: v(n-at) \geq e \} \\
    =  
    \begin{cases}
      p^{-\max\{0,e\}} \frac{T}{t} \phi(t)  \pm 2{d(t)} 
      & \text{if $p \nmid t$,} \\
      \frac Tt\phi(t)  \pm d(t)   & \text{if $p \mid t$ and $e\leq0$,} \\
      0                            & \text{if $p \mid t$ and $e>0$.} 
    \end{cases}
  \end{multline*}
\end{lemma}
\begin{proof}
  Clearly, if $e<0$, we may replace the condition $v(n-at)\geq e$ with
  $v(n-at)\geq0$; in other words, we may suppose $e$ non negative.

  Suppose $p\mid t$: if $e>0$, the statement is obvious; if $e=0$, it
  follows from Proposition~\ref{p:3}.
  Suppose now that $p \nmid t$. Then
  \begin{align}
    & \#\{ n\in\Z : 0\leq n \leq T, \; \gcd(n,t)=1, \; v(n-at) \geq e\} 
    \nonumber\\ 
    & \qquad = \# \left\{
      \begin{aligned}
        n = p^e n' + at :\; & n'\in\Z, \;  
        -\frac{at}{p^e} \leq n' \leq \frac{T-at}{p^e}, \; \\
        & \gcd(p^e n' + at,t)=1, \;
        v(p^e n') \geq e
      \end{aligned}
    \right\} 
    \nonumber \\
    & \qquad = \# \left\{
      n'\in\Z : -\frac{at}{p^e} \leq n' \leq \frac{T-at}{p^e}, \; 
      \gcd(n',t) = 1
    \right\}. 
    \label{eq:2}\\
    \intertext{Suppose $a\geq0$ and $T\geq at$, then 
      equation~\eqref{eq:2} becomes}
    & \qquad = \phi\left(t, \frac{T-at}{p^e}\right)
    + \phi\left(t, \frac{at}{p^e}\right);
    \nonumber \\
    \intertext{by proposition~\ref{p:3}, this is}
    & \qquad = p^{-e} \frac{T}{t} \phi(t) \pm 2{d(t)}.
    \label{eq:3}
  \end{align}
  If $T<at$ or $a<0$, it
  is a trivial calculation to verify that equation~\eqref{eq:3} still
  holds true.
\end{proof}
\begin{rem}
  The Lemma holds even if $a\in\Z_p$: it suffices to write $a = a' +
  \bigoh{p^\eta}$ with $a'\in\Z$ and $\eta$ large enough so that, for every
  positive $n\leq T$, $v(n-at) = v(n-a't)$. In particular, it holds if
  $a\in\Q$ with $v_p(a)\geq0$.
\end{rem}
\begin{lemma}
  \label{l:2}
  Suppose $m$ and $n$ are relatively prime integers. Then
  \begin{equation*}
    v(m/n) \geq e \qquad \text{if and only if} \qquad
    \begin{cases} 
      v(m) \geq  e & \text{if $e>0$,}\\
      v(n) \leq -e & \text{if $e\leq0$.}
    \end{cases}
  \end{equation*} 
\end{lemma}
\begin{proof}
  Obvious.
\end{proof}
\begin{prop}
  \label{lemma:1}
  For any $T>0$ we have
  \begin{align*}
    \sum_{n\leq T} \phi(n)
    &= \frac{1}{2\zeta(2)}T^2 + \bigoh{T\log T}, &
    \sum_{n\leq T} d(n) &= T\log T + \bigoh T, \\
    \sum_{\substack{
        n\leq T \\ p \mid n
      }} \phi(n) 
    &= \frac{1}{2\zeta(2)(p+1)} T^2 + \smalloh{T^2}, &
    \sum_{\substack{
        n\leq T \\ p \nmid n
      }} \phi(n) 
    &= \frac{p}{2\zeta(2)(p+1)} T^2 + \smalloh{T^2}.
  \end{align*}
\end{prop}
\begin{proof}
  The first and second assertions are very well known (see, e.g.,
  \cite[Chapter 3]{Apostol-IntroANT}) while the third clearly follows from
  the last one. Let
  \begin{equation*}
    \phi'(n) = 
    \begin{cases}
      \phi(n) & \text{if $p\nmid n$,} \\
      0 & \text{if $p\mid n$.}
    \end{cases}
  \end{equation*}
  Since $\phi'$, as well as $\phi$, is multiplicative, we have for $\Re(s)>2$:
  \begin{equation*}
    \sum^{\infty}_{\substack{
        n=1 \\ p \nmid n
      }} \frac{\phi(n)}{n^s} =
    \sum_{n=1}^{\infty} \frac{\phi'(n)}{n^s} =
    \frac{p^s-p}{p^s-1} \sum_{n=1}^\infty \frac{\phi(n)}{n^s} =
    \frac{p^s-p}{p^s-1} \cdot\frac{\zeta(s-1)}{\zeta(s)}.
  \end{equation*}
  In particular, the \textsc{lhs} is regular on the line $\Re(s)=2$ with
  the exception of a pole of first order at $s=2$ with residue
  $p/(p+1)\zeta(2)$. Thanks to the Tauberian theorem for Dirichlet series
  (see \cite[XV, \S3]{Lang-ANT}) it follows that, after a change of
  coordinate $s \mapsto s-1$,
  \begin{equation*}
    \sum_{n=1}^T \frac{\phi'(n)}{n} = \frac{p}{p+1}\cdot\frac{T}{\zeta(2)} +
    \smalloh{T}.
  \end{equation*}
  Applying the following ``integration'' lemma, we get the proposition.
  Note that, actually, the error term for the two last assertions is
  $\bigoh{T\log T}$: since we will not need this improved estimate, we
  will content ourselves with the simpler previous proof.
\end{proof}
\begin{lemma}
  \label{lemma:integration}
  Let $\{b_n\}_{n\in\N}$ be complex numbers and let $B(T) = \sum_{n=1}^T
  b_n$. Suppose that there is $\beta\in\C$ such that $B(T) = \beta T +
  \smalloh{T}$. Then
  \begin{equation*}
    \sum_{n=1}^T n b_n = \frac{\beta}{2} T^2 + \smalloh{T^2}.
  \end{equation*}
\end{lemma}

\section{Slices}
\label{sec:slices}

In order to prove Theorem~\ref{theor:1}, we count how many rational points in
$B(x,p^e)$ have a given height, then we use the previous Lemma to sum over
all heights.

\begin{df}
  For any $x\in\Q_p$ and $e\in\Z$, let
  \begin{math}
    B(x,p^{e};T) = B(x,p^e) \cap \Q(T),
  \end{math}
  where 
  \begin{math}
    \Q(T) = \{ t \in \Q : H(t) = T \}.
  \end{math}
\end{df}

\begin{lemma}
  \label{lemma:2}
  For every positive integer $T$: $\# \Q(T) = 4\phi(T)$.
\end{lemma}
\begin{proof}
  Obvious
\end{proof}

\begin{prop}
  \label{p:5}
  Fix $x\in\Q_p$ and $e\in \Z$ such that $v(x) \geq e$, then for any $t
  \in \Z^{>0}$ we have:
  \begin{enumerate}
  \item if $e>0$,
    \begin{equation*}
      \# B(x, p^{-e}; t) =
      \begin{cases}
        2p^{-e}\phi(t) \pm 4{d(t)} & \text{if $v(t)=0$,}\\
        0 & \text{if $0<v(t)<e$,}\\
        2\phi(t) & \text{if $v(t)\geq e$;}
      \end{cases}
    \end{equation*}
  \item if $e\leq0$,
    \begin{equation*}
      \# B(x, p^{-e}; t) =
      \begin{cases}
        2\bigl(2-p^{e-1}\bigr) \phi(t) \pm 4{d(t)} &
        \text{if $v(t)=0$,}\\ 
        4\phi(t) & \text{if $0<v(t) < 1-e$,}\\
        2\phi(t) & \text{if $v(t) \geq 1-e$.}
      \end{cases}
    \end{equation*}
  \end{enumerate}
\end{prop}
\begin{proof}
  Since $0\in B(x, p^{-e})$, we have $B(x, p^{-e}) = B(0, p^{-e})$ and
  thus, for any $t$,
  $B(x, p^{-e}; t) = B(0, p^{-e}; t)$. Therefore, we may as
  well suppose that $x=0$. Write $B(0, p^{-e}; t)$ as
  \begin{equation}
    \label{eq:4}
    \left\{
      \frac{m}{n} : (m,n) \in \Z\times\Z^{>0},\; \max\{ \abs{m}, n \} = t, \;
      \gcd(m,n)=1, \; v(m/n) \geq e
    \right\}.
  \end{equation}
  Suppose $e>0$: then~\eqref{eq:4} becomes, using Lemma~\ref{l:2},
  \begin{align*}
    & B(0,p^{-e}; t) \\
    & \qquad= \left\{
      m/n : (m,n) \in \Z\times\Z^{>0},\; \max\{ \abs{m}, n \} = t, \;
      \gcd(m,n)=1, \; v(m) \geq e
    \right\} \\
    & \qquad
    = \{±t/n : n\in\Z,\; 1\leq n \leq t,\; \gcd(n,t)=1,\; v(t) \geq e\} \\ 
    & \qquad\qquad \cup \; 
    \{m/t : m\in\Z,\; -t\leq m\leq t,\; \gcd(m,t)=1, \; v(m) \geq e\}.
  \end{align*}
  The first part of the proposition now follows from Lemma~\ref{l:1}.

  In order to prove the second part, notice that
  \begin{equation*}
    \Q(T) = B(0,p^{-e}; t) \cup \{r\in\Q : H(r)=T,\, v(r) \leq e-1 \}
  \end{equation*}
  and that $r\mapsto 1/r$ induces a 1-to-1 correspondence
  \begin{equation*}
      \{r\in\Q : H(r)=T,\, v(r) \leq e-1 \} \longleftrightarrow
      \{r\in\Q : H(r)=T,\, v(r) \geq 1-e \} 
  \end{equation*}
  Thus
  \begin{math}
    \# B(0,p^{-e}; t) = \#\Q(T) - \# B(0,p^{-(1-e)}; t)
  \end{math}
  and part (2) follows from part (1) and Lemma~\ref{lemma:2}.
\end{proof}

\begin{prop}
  \label{p:6}
  Fix $x\in\Q_p$ and $e\in \Z$ such that $v(x) < e$, then for
  any $t \in \Z^{>0}$ we have:
  \begin{enumerate}
   \item if $v(x) > 0$,
    \begin{equation*}
      \# B(x, p^{-e}; t) = 
      \begin{cases}
        2p^{-e}\phi(t) \pm 4{d(t)} & \text{if $v(t)=0$,}\\
        \displaystyle 2\frac{p^{1+v(x)-e}}{p-1} \phi(t) \pm 4{d(t)} 
        & \text{if $v(t)=v(x)$,} \\
        0 & \text{otherwise;}
      \end{cases}
    \end{equation*}
   \item if $v(x) = 0$,
    \begin{equation*}
      \#  B(x, p^{-e}; t) = 
      \begin{cases}
        4p^{-e} \phi(t) \pm 8{d(t)} & \text{if $v(t)=0$,}\\
        0 & \text{if $v(t) \neq 0$;}
      \end{cases}
    \end{equation*}
   \item if $v(x) < 0$,
    \begin{equation*}
      \#  B(x, p^{-e}; t) = 
      \begin{cases}
        2p^{2v(x)-e} \phi(t) \pm 4{d(t)} & \text{if $v(t)=0$,}\\
        \displaystyle 2\frac{p^{1+v(x)-e}}{p-1} \phi(t) \pm 4{d(t)} 
        & \text{if $v(t) = -v(x)$,}\\
        0 & \text{otherwise.}
      \end{cases}
    \end{equation*}
  \end{enumerate}
\end{prop}
\begin{proof}
  We have that
  \begin{align*}
    & B(x,p^{-e};t) \notag \\
    & \qquad = \left\{
      \begin{aligned}
        m/n :\; & (m,n) \in \Z\times\Z^{>0},\; 
        \max\{ \abs{m}, n \} = t, \; 
        \gcd(m,n)=1,\; \\
        & v\left(\frac{m}{n} - x \right) \geq e 
      \end{aligned}
    \right\}
    \notag \\[4pt]
    & \qquad = \left\{
      ±t/n :n\in\Z,\; 1\leq n \leq t,\; \gcd(n,t)=1,\; 
      v\left(\frac{±t}{n} - x\right) \geq e
    \right\} \\ 
    & \qquad\qquad \bigcup \;  \left\{
      m/t : m\in\Z,\; -t\leq m\leq t,\; \gcd(m,t)=1, \;
      v\left(\frac{m}{t} - x\right) \geq e
    \right\}. \notag
  \end{align*}
  Let us call these two sets respectively $B_1$ and $B_2$. 
  
  Consider $B_2$. We have
  \begin{equation}
    \label{eq:7}
    v\left(\frac{m}{t} - x\right) \geq e 
    \Longleftrightarrow
    v(m - xt) \geq e + v(t).
  \end{equation}

  Suppose $v(t)=0$. 
  \begin{itemize}
  \item If $v(x)\geq0$, we can apply Lemma~\ref{l:1}: since $e>v(x)$,
    we get
    \begin{math}
      \# B_2 = 2p^{-e} \phi(t) \pm 4d(t). 
    \end{math}
  \item If $v(x)<0$, we get
    \begin{math}
      v(m-xt) = v(x) < e,
    \end{math}
    hence Eq.~\eqref{eq:7} is false and $\# B_2 = 0$.
  \end{itemize}
  
  Suppose now $v(t)>0$. Then $\gcd(m,t)=1$ implies that $v(m)=0$.
  \begin{itemize}
  \item If $v(xt)>0$, then $v(m-xt) = v(m) = 0$ with $0 < v(x) + v(t) <
    e + v(t)$. Hence Eq.~\eqref{eq:7} is false.
  \item If $v(xt)=0$, let $\eta = v(t) > 0$, $x' = xp^\eta$ and $t' =
    tp^{-\eta}$. Hence
    \begin{equation*}
      \# B_2 = \# \{m : -t\leq m\leq t, \gcd(m,t)=1,
      v(m-x't')\geq e+\eta\}.
    \end{equation*}
    We can replace the condition $\gcd(m,p^\eta t')=1$ with
    $\gcd(m,t')=1$, since $p\mid m$ would imply $v(m-x't')= 0 = v(x) + v(t) <
    e + v(t)$. Lemma~\ref{l:1} yields
    \begin{equation*}
      \#B_2 = 2p^{-e-\eta} \frac{t}{t'} \phi(t') \pm 4d(t)
      = 2 \frac{p^{1-e-\eta}}{p-1} \phi(t) \pm 4d(t).
    \end{equation*}
    
  \item If $v(xt)<0$, then 
    \begin{math}
      v(m-xt) = v(xt) = v(x) + v(t) < e + v(t);
    \end{math}
    Eq.~\eqref{eq:7} is therefore false.
  \end{itemize}

  Putting everything together, we have
  \begin{equation*}
    \# B_2 = 
    \begin{cases}
      2 p^{-e} \phi(t) \pm 4{d(t)} &
      \text{if $v(t)=0$ and $v(x)\geq0$,}\\
      2 \frac{p^{1-e-v(t)}}{p-1} \phi(t) \pm 4d(t) &
      \text{if $v(t)>0$ and $v(x) = -v(t)$,}\\
      0 & \text{otherwise}.
    \end{cases}
  \end{equation*}
  
  \bigskip
    
  Consider now $B_1$. Write $x = x' p^\eta$ with $x'\in\Z_p$
  and $\eta=v(x)$.
  Since $v(x) < e$ and $v(\pm t/n - x) \geq e$ we have
  \begin{equation}
    \label{eq:6}
    v(t/n) = v(x) = \eta, \quad \text{with the constraint }
    \gcd(t,n) = 1
  \end{equation}

  Assume that $\eta \geq 0$, then Eq.~\eqref{eq:6} implies
  that $v(t)=\eta$ and $v(n)=0$; in particular $B_1 = \emptyset$ if
  $v(t) \neq v(x)$. Suppose thus $v(t)=\eta$ and write $t=t' p^\eta$. Then
  \begin{equation*}
    v\left(±\frac{t}{n} - x\right)
    = \eta + v\left(n x' \mp t'\right)
    = \eta +v(n \mp x'^{-1}t')
  \end{equation*}
  and, since $v(n)=0$,
  \begin{equation*}
    B_1 = \left\{
      ±t/n : n\in\Z,\; 1\leq n \leq t,\;
      \gcd(n,t')=1,\; v(n \mp x'^{-1}t') \geq e - \eta
    \right\}.
  \end{equation*}
  It follows, by lemma~\ref{l:1}, that
  \begin{equation*}
    \# B_1 
    = 2 p^{\eta-e}  \frac{t}{t'} \phi(t') \pm 4{d(t)}
    = \begin{cases}
      2 \frac{p^{1+\eta-e}}{p-1} \phi(t) \pm 4{d(t)}
      & \text{if $v(t)=v(x)>0$,}\\
      2 p^{-e} \phi(t) \pm 4{d(t)} & \text{if $v(t)=v(x)=0$,} \\
      0 & \text{otherwise.}
    \end{cases}
  \end{equation*}

  Assume now that $\eta < 0$, then equation~\eqref{eq:6}
  implies that $v(t) = 0$ and $v(n) = -\eta$; in particular, $B_1 =
  \emptyset$ if $v(t) > 0$. Suppose not and write $n = n'
  p^{-\eta}$. Then
  \begin{equation*}
    v\left(±\frac{t}{n} - x\right) = \eta + v(\pm t - n' x').
  \end{equation*}
  Since $e-\eta>0$ and $v(x'^{-1}t)=0$, we have
  \begin{equation*}
    \begin{split}
      \# B_1 & = \# \left\{
        \begin{aligned}
          ±t/n'p^{-\eta} : \; & n' \in \Z, \;
          1 \leq n' \leq p^\eta t, \;
          \gcd(n', t) = 1, \; \\
          & v(n' \mp x'^{-1}t ) \geq e - \eta
        \end{aligned}
      \right\}
      \\
      & = 2 p^{\eta-e} \frac{p^\eta t}{t} \phi(t) \pm 4{d(t)}
      = 2 p^{2\eta - e} \phi(t) \pm 4{d(t)}.
    \end{split}
  \end{equation*}
  The Proposition follows.
\end{proof}

\section{Proof of Theorem~\ref{theor:1}}
\label{sec:proof-theorem}

Let $e$ be a strictly positive integer. Then, by Proposition~\ref{lemma:1},
\begin{equation}
  \label{eq:5}
  \sum_{\substack{
      t\leq T \\ p^e \mid t
    }} \phi(t) =
  p^{e-1} \sum_{\substack{
      t\leq T/p^{e-1} \\ p \mid t
    }} \phi(t)
  = \frac{p^{1-e}}{2\zeta(2)(p+1)} T^2 + \bigoh{T\log T}.
\end{equation}
It follows that, the sum being over positive terms,
\begin{equation}
  \label{eq:9}
  \sum_{\substack{
      t\leq T \\ v(t) = e
    }} \phi(t) =
  \sum_{\substack{
      t\leq T \\ p^e \mid t
    }} \phi(t)
  - \sum_{\substack{
      t\leq T \\ p^{e+1} \mid t
    }} \phi(t)
  = p^{-e} \frac{p-1}{2\zeta(2)(p+1)} T^2 + \bigoh{T\log T}.
\end{equation}
Suppose now that $0<e\leq v(x)$. Then Proposition~\ref{p:5} yields
\begin{align*}
  \sum_{t\leq T} \# B(x,p^{-e}; t) & =
  \sum_{\substack{
      t\leq T \\ p \nmid t
    }} \left(
    \frac{2}{p^e}\phi(t) \pm 4d(t)
  \right)
  + \sum_{\substack{
      t\leq T \\ p^e \mid t
    }} 2\phi(t)
  \\
  &= \frac{2p^{1-e}}{2\zeta(2)(p+1)} T^2 +
  \frac{2p^{1-e}}{2\zeta(2)(p+1)} T^2 + \bigoh{T\log T}
  \\
  &= \frac{4p^{1-e}}{2\zeta(2)(p+1)} T^2 + \bigoh{T\log T};
\end{align*}
therefore
\begin{equation*}
  \mu\bigl(B(x,p^{-e})\bigr) =
  \lim_{T\to+\infty} \frac{
    \sum_{t\leq T} \# B(x,p^{-e}; t)}{
    \sum_{t\leq T} \# \Q(t)}
  = \frac{p^{1-e}}{p+1}.
\end{equation*}
If $e\leq0$ and $e\leq v(x)$ we have, instead:
\begin{multline*}
  \sum_{t\leq T} \# B(x,p^{-e};t) =
  \sum_{t\leq T} 4\phi(t)
  -2 \sum_{\substack{ t\leq T \\ p \nmid t}}\left(
    p^{e-1}\phi(t) \pm 4d(t)
  \right)
  -2 \sum_{\substack{ t\leq T \\ p^{1-e} \mid t}} \phi(t) \\
  = \frac{4}{2\zeta(2)}T^2 -\frac{2p^e}{2\zeta(2)(p+1)}T^2
  - \frac{2p^e}{2\zeta(2)(p+1)}T^2 + \bigoh{T\log T} \\
  = 4\left(1 - \frac{p^e}{p+1}\right) \frac{1}{2\zeta(2)}T^2 +
  \bigoh{T\log T}.
\end{multline*}
Thus 
\begin{equation*}
  \mu\bigl(B(x,p^{-e})\bigr) = 1 - \frac{p^e}{p+1}.
\end{equation*}
Suppose now $e>v(x)>0$. Then Proposition~\ref{p:6} and Eq.~\eqref{eq:9}
yield
\begin{align*}
  \sum_{t\leq T} \# B(x,p^{-e}; t) & =
  \sum_{\substack{
      t\leq T \\ p \nmid t
    }} \left(
    \frac{2}{p^e}\phi(t) \pm 4d(t)
  \right)
  + \sum_{\substack{
      t\leq T \\ v(t) = v(x)
    }} \left(
    2\frac{p^{v(x)-e}}{1-1/p}\phi(t) \pm 4d(t)
  \right)
  \\
  &= \frac{4p^{1-e}}{2\zeta(2)(p+1)} T^2 + \bigoh{T\log T}; \\
\end{align*}
so that
\begin{math}
  \mu\bigl( B(x, p^{-e})\bigr) = p^{1-e}/(p+1).
\end{math}
If $e>v(x)=0$, the calculation clearly gives the same result.
Suppose at last that $e>v(x)$ with $v(x)<0$. Then
\begin{align*}
  \sum_{t\leq T} \# B(x,p^{-e}; t) & =
  \sum_{\substack{
      t\leq T \\ p \nmid t
    }} \left(
    2p^{2v(x)-e} \phi(t) \pm 4d(t)
  \right)
  \\
  & \quad + \sum_{\substack{
      t\leq T \\ v(t)=-v(x)
    }} \left(
    2 \frac{p^{1+v(x)-e}}{p-1} \phi(t) \pm 4d(t)
  \right)
  \\
  &= \frac{4p^{1+2v(x)-e}}{2\zeta(2)(p+1)} T^2 + \bigoh{T\log T};
  \phantom{\frac{\Biggl(}{}}
\end{align*}
hence 
\begin{math}
  \mu\bigl( B(x, p^{-e})\bigr) = p^{1+2v(x)-e}/(p+1).
\end{math}
\qed

\bibliography{biblio}
\bibliographystyle{abbrv}
\end{document}